\documentclass[12pt]{amsart}
\usepackage{amsmath, amssymb, latexsym, amsthm}
\usepackage{mathrsfs}

\newcommand{\Bc}[9]{\bibitem{#1} {#2}, \emph{#3}, in: \textbf{#4} (#5), #6 #7, #8--#9.}

\newcommand{\ed}{

\end{document}
}

\newcommand{\next}[1]{{#1^+}}
\newcommand{\Menger}{\mathrm{Menger}}
\newcommand{\intvl}[2]{{[#1(#2),#1(#2\!+\!1))}}

\newcommand{\bq}{\begin{quote}}
\newcommand{\eq}{\end{quote}}

\newcommand{\Cantor}{{\{0,1\}^\N}}

\newcommand{\N}{\mathbb{N}}
\newcommand{\NN}{{\N^{\N}}}

\newcommand{\PN}{{P(\N)}}
\newcommand{\roth}{{[\N]^{\aleph_0}}}

\newcommand{\seq}[1]{\{#1\}_{n\in\N}}

\newcommand{\cF}{\mathcal{F}}

\newcommand{\cU}{\mathcal{U}}
\newcommand{\Union}{\bigcup}

\newcommand{\Impl}{\Rightarrow}
\long\def\forget#1\forgotten{}

\newcommand{\fb}{\mathfrak{b}}

\newcommand{\fd}{\mathfrak{d}}
\newcommand{\fg}{\mathfrak{g}}

\newcommand{\x}{\times}
\newcommand{\Iff}{\Leftrightarrow}
\newcommand\comp{^{\text{\tt c}}}
\newcommand{\nin}{\not\in}

\newcommand{\sbst}{\subseteq}

\newcommand{\sm}{\setminus}
\newcommand{\as}{\subseteq^*}

\newcommand{\add}{\mathsf{add}}

\newcommand{\cf}{\mathsf{cf}}

\newtheorem{thm}{Theorem}
\newtheorem{prop}[thm]{Proposition}

\newtheorem{lem}[thm]{Lemma}
\newtheorem{cor}[thm]{Corollary}

\theoremstyle{definition}

\theoremstyle{remark}

\newcommand{\be}{\begin{enumerate}}
\newcommand{\ee}{\end{enumerate}}
\newcommand{\bi}{\begin{itemize}}
\newcommand{\ei}{\end{itemize}}

\forget
\forgotten

\author{Boaz Tsaban}
\thanks{Supported by the Koshland Center for Basic Research.}
\address{Boaz Tsaban, Department of Mathematics,
Weizmann Institute of Science,
Rehovot 76100,
Israel}
\email{boaz.tsaban@weizmann.ac.il}
\urladdr{http://www.cs.biu.ac.il/\~{}tsaban}

\author{Lyubomyr Zdomskyy}
\address{Lyubomyr Zdomskyy, Department of Mechanics and Mathematics,
Iv\-an Franko Lviv National University,
Universytetska 1, Lviv 79000, Ukraine.}
\email{lzdomsky@rambler.ru}

\title[Menger and groupwise density]{Menger's covering property and groupwise density}

\begin{document}

\begin{abstract}
We establish a surprising connection between Men\-ger's classical
covering property and Blass-Laflamme's modern combinatorial notion of groupwise density.
This connection implies a short proof of the groupwise density bound
on the additivity number for Menger's property.
\end{abstract}


\maketitle

\section{Introduction and basic facts}

Unless otherwise indicated, all spaces considered in this paper
are assumed to be separable, zero-dimensional, and metrizable.
Consequently, all open covers may be assumed to be countable.

\emph{Menger's property} (1924), defined in \cite{Menger24},
is a generalization of $\sigma$-compactness.
The following more familiar reformulation of this property was given
by Hurewicz in \cite{Hure25}:
A space $X$ has Menger's property if, and only if, for each sequence
$\seq{\cU_n}$ of open covers of $X$, there exist finite sets
$\cF_n\sbst\cU_n$, $n\in\N$, such that $\Union_{n\in\N} \cF_n$ is a cover of $X$.

Hurewicz's reformulation easily implies that Menger's property
is preserved under continuous images and is hereditary for closed subsets.
It is also not difficult to see that it is preserved under
taking countable unions.

\begin{cor}\label{Fsigmahered}
Menger's property is hereditary for $F_\sigma$ subsets.\hfill\qed
\end{cor}

The following preservation property will also be useful
(see \cite{AddQuad} for a proof).

\begin{prop}[folklore]\label{Ksigma}
Assume that $X$ has Menger's property and $K$ is $\sigma$-compact.
Then $X\x K$ has Menger's property.\hfill\qed
\end{prop}

In 1927, Hurewicz has essentially obtained the following combinatorial
characterization of Menger's property.
View $\N$ as a discrete topological space,
and endow the \emph{Baire space} $\NN$ with the Tychonoff product topology.
Define a partial order $\le^*$ on $\NN$ by:
$$f\le^* g\quad \mbox{if}\quad f(n)\le g(n)\mbox{ for all but finitely many }n.$$
A subset $D$ of $\NN$ is \emph{dominating} if
for each $g\in\NN$ there exists $f\in D$ such that $g\le^* f$.

\begin{thm}[Hurewicz \cite{Hure27}]\label{hure}
A space $X$ has Menger's property if, and only if, no
continuous image of $X$ in $\NN$ is dominating.\hfill\qed
\end{thm}

While the proof of this assertion is very easy \cite{Rec94},
this characterization has found numerous applications
(see \cite{LecceSurvey, ict} and references therein).

An important application of the Hurewicz Theorem \ref{hure} is the following.
Let $\add(\Menger)$ denote the minimal cardinality of a family
of spaces with Menger's property, whose union does not have Menger's property.
Let $\fb$ denote the minimal cardinality of
an unbounded (with respect to $\le^*$) subset of $\NN$,
and $\fd$ denote the minimal cardinality of
a dominating subset of $\NN$.
By Theorem \ref{hure}, the minimal cardinality of a space which
does not have Menger's property is $\fd$.
Using this and Theorem \ref{hure} again, we have that
$\fb\le\add(\Menger)\le\cf(\fd)$.

In this paper we give a new characterization of Menger's property,
in terms of a combinatorial property whose connection to Menger's property is
less transparent.
To this end, we need several more definitions.

Let $\roth$ denote the collection of all infinite sets of natural numbers.
For $a,b\in\roth$, $a$ is an \emph{almost subset} of $b$, $a\as b$, if
$a\sm b$ is finite. A family $G\sbst\roth$ is \emph{groupwise dense} if
it contains all almost subsets of its elements, and for each
partition of $\N$ into finite intervals, there is an infinite set of intervals
in this partition whose union is a member of $G$.

Intuitively, groupwise dense families are large.
Roughly speaking, our main result asserts that
if a space $X$ has Menger's property, then for each continuous image of
$X$ in $\NN$ there are ``groupwise-densely'' many functions witnessing
that it is not dominating.

The \emph{groupwise density number} $\fg$ is the minimal cardinality of a
collection of groupwise dense families whose intersection is empty.
This relatively new notion is due to Blass and Laflamme \cite{BlLaf89, Blass90}.
It follows at once that $\fg\le\add(\Menger)$.
This consequence,
which was previously obtained using much more complicated arguments \cite{SF1},
strengthens the Blass-Mildenberger result that $\fg\le\cf(\fd)$ \cite{BlassMil99}.

\section{A new characterization of Menger's property}\label{main}

$\roth$ is a topological subspace of $\PN$, where the topology on
$\PN$ is defined by identifying it with the Cantor space $\Cantor$.
For $m<n$, define $[m,n)=\{m,m+1,\dots,n-1\}$.
For each finite $F\sbst\N$ and each $n\in\N$, define
$$O_{F,n}=\{a\in\PN : a\cap[0,n)=F\}.$$
The sets $O_{F,n}$ form a clopen basis for the topology on $\PN$.

For $a\in\roth$, define an element $\next{a}$ of $\NN$
by
$$\next{a}(n) = \min\{k\in a : n<k\}$$
for each $n$. The following is the main result of this paper.

\begin{thm}\label{gd}
A space $X$ has Menger's property if, and only if,
for each continuous image $Y$ of $X$ in $\NN$,
the family
$$G=\{a\in\roth : (\forall f\in Y)\ \next{a}\not\le^* f\}$$
is groupwise dense.
\end{thm}
\begin{proof}
The direction $(\Leftarrow)$ is immediate from the Hurewicz Theorem \ref{hure}.
We prove the more interesting direction $(\Impl)$.

Assume that $Y$ is a continuous image of $X$ in $\NN$.
Then $Y$ has Menger's property.
By Proposition \ref{Ksigma}, $\PN\x Y$ has Menger's property.

\begin{lem}\label{Clemma}
The set
$$C = \{(a,f)\in\roth\x\NN : \next{a}\le^* f\}$$
is an $F_\sigma$ subset of $\PN\x\NN$.
\end{lem}
\begin{proof}
Note that
$$C = \Union_{m\in\N}\bigcap_{n\ge m}\{(a,f)\in \PN\x\NN : (n,f(n)]\cap a\neq\emptyset\}.$$
(The nonempty intersection for infinitely many $n$ allows the replacement of $\roth$ by $\PN$.)
For fixed $m$ and $n$, the set $\{(a,f)\in P(\N)\x\NN : (n,f(n)]\cap a\neq\emptyset\}$ is clopen:
Indeed, if $\lim_k(a_k,f_k)=(a,f)$ then for all large enough $k$, $f_k(n)=f(n)$,
and therefore for all larger enough $k$, $(n,f_k(n)]\cap a_k=(n,f(n)]\cap a$.
Thus, $(a_k,f_k)$ is in the set if, and only if, $(a,f)$ is in the set.
\end{proof}

By Corollary \ref{Fsigmahered} and Lemma \ref{Clemma},
$C\cap (\PN\x Y)$ has Menger's property, and
therefore so does its projection $Z$ on the first coordinate.
By the definition of $Z$, $G=Z\comp$, the complement of $Z$ in $\roth$.
Note that $G$ contains all almost subsets of its elements.

For $a\in\roth$ and an increasing $h\in\NN$, define
$$a/h = \{n : a\cap \intvl{h}{n}\neq\emptyset\}.$$
For $S\sbst\roth$, define $S/h = \{a/h : a\in S\}$.

\begin{lem}\label{tala}
Assume that $G\sbst\roth$ contains all almost subsets of its elements.
Then: $G$ is groupwise dense if, and only if, for each increasing $h\in\NN$,
$G\comp/h\neq\roth$.
\end{lem}
\begin{proof}
For each increasing $h\in\NN$ and each $a\in\roth$,
$$\Union_{n\in a}\intvl{h}{n}\nin G\Iff \Union_{n\in a}\intvl{h}{n}\in G\comp
\Iff a\in G\comp/h.$$
The lemma follows directly from that.
\end{proof}

Assume that $G$ is not groupwise dense.
By Lemma \ref{tala}, there is an increasing $h\in\NN$ such that
$Z/h = G\comp/h=\roth$.
The natural mapping $\Psi:Z\to Z/h$ defined by $\Psi(a)=a/h$ is
a continuous surjection. It follows that $\roth$ has Menger's property.
But this is absurd:
The image of $\roth$ in $\NN$, under the
continuous mapping assigning to each $a\in\roth$ its increasing enumeration,
is a dominating subset of $\NN$.
Thus, $\roth$ does not have Menger's property -- a contradiction.
\end{proof}

We immediately obtain the following.

\begin{cor}[\cite{SF1}]
Each union of less than $\fg$ many spaces having Menger's property,
has Menger's property.
\end{cor}
\begin{proof}
Assume that $\kappa<\fg$ and for each $\alpha<\kappa$, $X_\alpha$ has Menger's property,
and that $X=\Union_{\alpha<\kappa}X_\alpha$.
By the Hurewicz Theorem \ref{hure}, it suffices to show that no continuous image of $X$ in
$\NN$ is dominating.
Indeed, assume that $\Psi:X\to\NN$ is continuous.
By Theorem \ref{gd}, for each $\alpha$ the family
$$G_\alpha =\{a\in\roth : (\forall f\in \Psi[X_\alpha])\ \next{a}\not\le^* f\}$$
is groupwise dense.
Thus, there exists $a\in\bigcap_{\alpha<\kappa}G_\alpha$.
Then $\next{a}$ witnesses that $\Psi[X]$ is not dominating.
\end{proof}

\section{Additional remarks}

The proof of Theorem \ref{gd} shows that if
a space $X$ has Menger's property, then
for each continuous image $Y$ of $X$ in $\NN$,
the family
$$G=\{a\in\roth : (\forall f\in Y)\ \next{a}\not\le^* f\}$$
is \emph{coMenger}, i.e., its complement in $\roth$ has Menger's property.
The proof also shows that coMenger sets in $\roth$ containing
all almost subsets of their elements are groupwise dense.

It is well known \cite{BlassHBK} that if $G\sbst\roth$
contains all almost subsets of its elements, then
$G$ is groupwise dense if, and only if, $G$ is nonmeager in
$\roth$. Thus, in Theorem \ref{gd}, ``groupwise dense'' can be
replaced by ``nonmeager''.

Using arguments similar to those in the proof of Theorem \ref{gd},
we have the following.

\begin{thm}
A space $X$ has Menger's property if, and only if,
for each continuous image $Y$ of $X$ in $\NN$,
the family
$$G=\{g\in\NN : (\forall f\in Y)\ g\not\le^* f\}$$
is nonmeager.\hfill\qed
\end{thm}

This is a structural extension of the same assertion
for spaces $X$ of cardinality smaller than $\fd$,
which was proved in \cite{SFH}.

\ed